\documentclass[11pt]{amsart}
\usepackage{amssymb}

\usepackage{graphicx}
\usepackage{amscd}
\usepackage{amsmath}

\newtheorem{theorem}{Theorem}
\theoremstyle{plain}
\newtheorem{acknowledgement}{Acknowledgement}

\newtheorem{corollary}{Corollary}

\newtheorem{definition}{Definition}
\newtheorem{example}{Example}

\newtheorem{lemma}{Lemma}

\newtheorem{proposition}{Proposition}
\newtheorem{remark}{Remark}

\numberwithin{equation}{section}

\begin{document}
\title[Absolutely Summing Mappings]{Cotype and non linear absolutely summing mappings}
\author{Daniel Pellegrino}
\address[Daniel Pellegrino]{Depto de Matem\'{a}tica e Estat\'{i}stica- UFPB Campus II- Caixa Postal
10044- Campina Grande-PB-Brazil and IMECC-UNICAMP- Caixa Postal 6065-
Campinas SP- Brazil}
\email[Daniel Pellegrino]{dmp@dme.ufpb.br}

\begin{abstract}
In this paper we study absolutely summing mappings on Banach spaces by
exploring the cotype of their domains and ranges. It is proved that every $n$%
-linear mapping from $\mathcal{L}_{\infty}$-spaces into $\mathbb{K}$ is $%
(2;2,...,2,\infty)$-summing and also shown that every $n$-linear mapping
from $\mathcal{L}_{\infty}$-spaces into $F$ is $(q;2,...,2)$-summing
whenever $F$ has cotype $q.$ We also give new examples of analytic summing
mappings and polynomial and multilinear versions of a linear Extrapolation
Theorem.
\end{abstract}

\maketitle

\section{Introduction}

In the fifties, A. Grothendieck%
%TCIMACRO{\UNICODE{0xb4}}%
%BeginExpansion
\'{}%
%EndExpansion
s seminal paper \cite{Grothendieck} ''Resum\'{e} de la th\'{e}orie
m\'{e}trique des produits tensoriels topologiques'' provided the
fundamentals of the absolutely summing operators theory. Subsequently, J.
Lindenstrauss and A. Pe\l czy\'{n}ski \cite{Lindenstrauss} simplified
Grothendieck%
%TCIMACRO{\UNICODE{0xb4}}%
%BeginExpansion
\'{}%
%EndExpansion
s tensorial notations leading to many interesting results. The multilinear
theory of absolutely summing mappings was outlined by Pietsch \cite{Pietsch2}
and has been developed by several authors (Alencar and Matos \cite
{AlencarMatos}, Floret and Matos \cite{MatosFloret}, Matos \cite{Matos2},
Schneider \cite{Schneider}, Tonge and Melendez \cite{Tonge}, Botelho \cite
{Botelho},\cite{Botelho1}, among others). Matos \cite{Matos2},\cite{Matos},
\cite{Matos1} also begun to study the concept of holomorphic absolutely
summing mappings and a more general definition in such a way that the origin
was not a distinguished point. The contribution of the notion of cotype to
this theory is relevant and can be seen in \cite{Botelho},\cite{Botelho1}
and \cite{MatosFloret}. In this paper, we will generalize several results of
\cite{Botelho1} and \cite{Botelho} and also give new Coincidence Theorems
and examples of absolutely summing holomorphic and analytic mappings.

\section{Notation, general concepts and basic results}

Throughout this paper $E,E_{1},...,E_{n},F,X,Y$ will always denote Banach
spaces and the scalar field $\mathbb{K}$ can be either $\mathbb{R}$ or $%
\mathbb{C}$. We will denote by $C(K)$ the Banach space of continuous scalar
valued functions on $K($compact Hausdorff space$)$ endowed with the $\sup$
norm.

The Banach space of all $n$-linear continuous mappings from $%
E_{1}\times...\times E_{n}$ into $F$ endowed with the canonical norm will be
denoted by $\mathcal{L}(E_{1},...,E_{n}$;$F)$ and the Banach space of all
continuous $n$-homogeneous polynomials $P$ from $E$ into $F$ with the norm $%
\Vert P\Vert=\sup\{\Vert Px\Vert;\Vert x\Vert\leq1\}$ will be denoted by $%
\mathcal{P}(^{n}E,F).$ A mapping $f:E\rightarrow F$ will be said analytic at
the point $a\in E,$ if there exist a ball $B_{\delta}(a)$ and a sequence of
polynomials $P_{k}\in\mathcal{P}(^{k}E,F)$ $\ $such that
\begin{equation*}
f(x)=\sum_{k=0}^{\infty}P_{k}(x-a)\text{ uniformly for }x\in B_{\delta}(a).
\end{equation*}

Henceforth $\delta_{a}$ will be called the radius of convergence of $f$
around $a.$ To emphasize the case $\mathbb{K}=\mathbb{C},$ we will sometimes
use the term ``holomorphic'' in the place of ``analytic''. Every analytic
mapping in the whole space will be called entire mapping.

For the natural isometry
\begin{equation*}
\Psi:\mathcal{L}(E_{1},...,E_{n};F)\rightarrow\text{ }\mathcal{L}%
(E_{1},...,E_{t};\mathcal{L}(E_{t+1},...,E_{n};F))
\end{equation*}
we will use the following convention: If $T\in\mathcal{L}(E_{1},...,E_{n};F)
$ then $\Psi(T)=T_{1}$ and if $T\in\mathcal{L}(E_{1},...,E_{t};\mathcal{L}%
(E_{t+1},...,E_{n};F)),$ then $\Psi^{-1}(T)=T_{0}.$

For $p\in]0,\infty\lbrack,$ the linear space of all sequences $%
(x_{j})_{j=1}^{\infty}$ in $E$ such that
\begin{equation*}
\Vert(x_{j})_{j=1}^{\infty}\Vert_{p}=(\sum_{j=1}^{\infty}\Vert x_{j}\Vert
^{p})^{\frac{1}{p}}<\infty
\end{equation*}
will be denoted by $l_{p}(E).$ We will also denote by $l_{p}^{w}(E)$ the
linear space of the sequences $(x_{j})_{j=1}^{\infty}$ in $E$ such that
\begin{equation*}
(<\varphi,x_{j}>)_{j=1}^{\infty}\in l_{p}(\mathbb{K}),
\end{equation*}
for every continuous linear functional $\varphi:E\rightarrow\mathbb{K}.$ We
also define $\Vert.\Vert_{w,p}$ in $l_{p}^{w}(E)$ by
\begin{equation*}
\Vert(x_{j})_{j\in\mathbb{N}}\Vert_{w,p}=\sup_{\varphi\in B_{E}%
%TCIMACRO{\UNICODE{0xb4}}%
%BeginExpansion
{\acute{}}%
%EndExpansion
}(\sum_{j=1}^{\infty}\mid<\varphi,x_{j}>\mid^{p})^{\frac{1}{p}}.
\end{equation*}
The case $p=\infty$ is just the case of bounded sequences and in $l_{\infty
}(E)$ we use the $\sup$ norm. The linear subspace of $l_{p}^{w}(E)$ of all
sequences $(x_{j})_{j=1}^{\infty}\in l_{p}^{w}(E),$ such that
\begin{equation*}
\lim_{m\rightarrow\infty}\Vert(x_{j})_{j=m}^{\infty}\Vert_{w,p}=0,
\end{equation*}
is a closed linear subspace of $l_{p}^{w}(E)$ and will be denoted by $%
l_{p}^{u}(E)$. The case $p=1$ motivates the name unconditionally $p$%
-summable sequences for the elements of $l_{p}^{u}(E)$ (\cite{Matos2}). One
can see that $\Vert.\Vert_{p}$ $(\Vert.\Vert_{w,p})$ is a $p$-norm in $%
l_{p}(E)($ $l_{p}^{w}(E))$ for $p<1$ and a norm in $l_{p}(E)($ $%
l_{p}^{w}(E)) $ for $p\geq1.$ In any case, they are complete metrizable
linear spaces.

\begin{definition}
Let $2\leq q\leq\infty$ and $(r_{j})_{j=1}^{\infty}$ be the Rademacher
functions. The Banach space $E$ has cotype $q,$ if there exists $%
C_{q}(E)\geq0,$ such that, for every $k\in\mathbb{N}$ and $%
x_{1},...,x_{k}\in E,$%
\begin{equation*}
(\sum_{j=1}^{k}\Vert x_{j}\Vert^{q})^{\frac{1}{q}}\leq C_{q}(E)(\int
\limits_{0}^{1}\Vert\sum_{j=1}^{k}r_{j}(t)x_{j}\Vert^{2}dt)^{\frac{1}{2}}.
\end{equation*}
To cover the case $q=\infty$, we replace $(\sum_{j=1}^{k}\Vert x_{j}\Vert
^{q})^{\frac{1}{q}}$ by
\begin{equation*}
\max\{\Vert x_{j}\Vert;1\leq j\leq k\}.
\end{equation*}
We will define the cotype of $E$ by $\cot E=\inf\{2\leq q\leq\infty;E$ has
cotype $q\}.$
\end{definition}

\begin{definition}
\label{deff}(Matos) A continuous $n$-linear mapping $T:E_{1}\times...\times
E_{n}\rightarrow F$ is absolutely $(p;q_{1},...,q_{n})$-summing (or $%
(p;q_{1},...,q_{n})$-summing) at $(a_{1},...,a_{n})\in E_{1}\times...\times
E_{n}$ if
\begin{equation*}
(T(a_{1}+x_{j}^{(1)},...,a_{n}+x_{j}^{(n)})-T(a_{1},...,a_{n}))_{j=1}^{%
\infty }\in l_{p}(F)
\end{equation*}
for every $(x_{j}^{(s)})_{j=1}^{\infty}\in l_{q_{s}}^{w}(E)$, $s=1,...,n.$ A
continuous $n$-homogeneous polynomial $P:E\rightarrow F$ is absolutely $(p;q)
$-summing (or $(p;q)$-summing) at $a\in E$ if
\begin{equation*}
(P(a+x_{j})-P(a))_{j=1}^{\infty}\in l_{p}(F)
\end{equation*}
for every $(x_{j})_{j=1}^{\infty}\in l_{q}^{w}(E).$
\end{definition}

The space of all $n$-homogeneous polynomials $P:E\rightarrow F$ which are $%
(p;q)$-summing (at every point) will be denoted by $\mathcal{P}%
_{as(p;q)(E)}(^{n}E;F).$ The space of all$\ n$-homogeneous polynomials $%
P:E\rightarrow F$ which are $(p;q)$-summing (at the origin) will be denoted
by $\mathcal{P}_{as(p;q)}(^{n}E;F).$ Analogously for $n$-linear mappings.

It must be noticed that the aforementioned definition, where the origin is
not a privileged point, is actually a more restrictive definition. For
example, if $n>1$ every $n$-linear mapping $T$ from $l_{1}\times...\times
l_{1}$ into $l_{1}$ is absolutely $(1;1)$-summing at the origin, but we can
always find $a\neq0$ such that $T$ is not absolutely $(1;1)$-summing at $a$
\cite{Matos1}$.$ Besides, the above definition turns possible to consider an
absolutely summing holomorphy type in the sense of Nachbin (see \cite{Matos}%
).

One can prove that if $r<s$ then the unique polynomial which is absolutely $%
(r;s)$-summing at every point is the trivial.

For $n$-homogeneous polynomials and $n$-linear mappings,\ the polynomials ($%
n $-linear mappings) $(\frac{p}{n};p)$-summing will be called $p$-dominated
polynomials\ ($n$-linear mappings) (see \cite{Matos2},\cite{Tonge}).\ For
the $p$-dominated polynomials ($n$-linear mappings) several natural versions
of linear results still hold, as well as Factorization Theorems, Domination
Theorem, etc. \cite{Matos2},\cite{Tonge},\cite{Schneider}.

The following characterization will be useful:

\begin{theorem}
(Matos \cite{Matos2}) Let $P$ be an $m$-homogeneous polynomial from $E$ into
$F$.Then, the following statements are equivalent:

(1) $P$ is absolutely $(p;q)$-summing at $0$.

(2)There exists $L>0$ such that
\begin{equation*}
(\sum_{j=1}^{\infty}\Vert P(x_{j})\Vert^{p})^{\frac{1}{p}}\leq L\Vert
(x_{j})_{j=1}^{\infty}\Vert_{w,q}^{m}\forall(x_{j})_{j=1}^{\infty}\in
l_{q}^{u}(E).
\end{equation*}

(3)There exists $L>0$ such that
\begin{equation*}
(\sum_{j=1}^{k}\Vert P(x_{j})\Vert^{p})^{\frac{1}{p}}\leq
L\Vert(x_{j})_{j=1}^{k}\Vert_{w,q}^{m}\forall k\in\mathbf{N,\forall}
x_{1},...,x_{k}.
\end{equation*}

(4) $(P(x_{j}))_{j=1}^{\infty}\in l_{p}(F)$ for every $(x_{j})_{j=1}^{\infty
}\in l_{q}^{w}(E).$
\end{theorem}

The infimum of the possible constants $L>0$ is a norm for the case $p\geq1$
or a $p$-norm for the case $p<1$(\cite{Matos2} or \cite{Pietsch} page 91) on
the space of the absolutely $(p;q)$-summing polynomials. In any case, we
have complete topological metrizable spaces. We will use the notation $\Vert
.\Vert_{as(p;q)}$ for this norm ($p-$norm).

The characterization for the multilinear case and the definition of the norm
($p$-norm) follows the same reasoning.

The following Theorem plays an important role in our future results:

\begin{theorem}
\label{Shera}(Maurey-Talagrand \cite{Talagrand})If $E$ has cotype $p,$ then $%
id:E\rightarrow E$ is $(p;1)$-summing. The converse is true, except for $p=2$%
.
\end{theorem}

The next definition, due to Lindenstrauss and Pe\l czy\'{n}ski is of
fundamental importance in the local study of Banach spaces and their
properties:

\begin{definition}
Let $1\leq p\leq\infty$ and let $\lambda>1.$ The Banach space $X$ is said to
be an $\mathcal{L}_{p,\lambda}$ space if every finite dimensional subspace $E
$ of $X$ is contained in a finite dimensional subspace $F$ of $X$ for which
there exists an isomorphism $v_{E}:F\rightarrow l_{p}^{\dim F}$ with $\Vert
v_{E}\Vert\Vert v_{E}^{-1}\Vert<\lambda.$ We say that $X$ is an $\mathcal{L}%
_{p}$ space if it\ is an $\mathcal{L}_{p,\lambda}$ space for some $\lambda>1.
$
\end{definition}

\section{Absolutely summing polynomials and multilinear mappings explored by
the cotype of their ranges}

The relation between cotype and absolutely summing linear mappings is clear
by Theorem \ref{Shera}. For points different from the origin we have the
straightforward following results:

\begin{lemma}
\label{Defant}Every continuous $n$-linear mapping $T:E_{1}\times...\times
E_{n}\rightarrow F$ is such that
\begin{equation*}
(T(a_{1}+x_{j}^{(1)},...,a_{n}+x_{j}^{(n)})-T(a_{1},...,a_{n}))_{j=1}^{%
\infty }\in l_{1}^{w}(F)
\end{equation*}
whenever $\ (x_{j}^{(1)})_{j=1}^{\infty}\in
l_{1}^{w}(E_{1}),...,(x_{j}^{(n)})_{j=1}^{\infty}\in l_{1}^{w}(E_{n}).$ The
polynomial version is immediate.
\end{lemma}

Proof. We just need to invoke a well known, albeit unpublished, result of
Defant and Voigt which states that every scalar valued $n$-linear mapping is
absolutely $(1;1)$-summing at the origin (see \cite{Matos}, Theorem 1.6 or
\cite{Matos2}), and explore multilinearity.

\begin{theorem}
\label{bbv}If $F$ has cotype $q$, then every continuous $n$-linear mapping
from $E_{1}\times...\times E_{n}$ into $F$ is $(q;1)$-summing on $%
E_{1}\times...\times E_{n}$. The polynomial case is also valid.
\end{theorem}

Proof. Since $F$ has cotype $q$, Theorem \ref{Shera} and Lemma \ref{Defant}
provide
\begin{equation*}
(\sum_{j=1}^{\infty}\Vert
T(a_{1}+x_{j}^{1},...,a_{n}+x_{j}^{n})-T(a_{1},...,a_{n})\Vert^{q})^{\frac{1%
}{q}}\leq
\end{equation*}
\begin{equation*}
\leq\Vert(T(a_{1}+x_{j}^{1},...,a_{n}+x_{j}^{n})-T(a_{1},...,a_{n}))_{j=1}^{%
\infty}\Vert_{w,1}<\infty
\end{equation*}
whenever $(x_{j}^{1})_{j=1}^{\infty}\in
l_{1}^{w}(E_{1}),...,(x_{j}^{n})_{j=1}^{\infty}\in
l_{1}^{w}(E_{n}).\blacksquare$

Theorem \ref{bbv} generalizes to points different from the origin the
following result:

\begin{theorem}
(Botelho \cite{Botelho}) If $F$ has cotype $q$ then every continuous $n$%
-homogeneous polynomial from $E$ into $F$ is $(q;1)$-summing at the origin.
\end{theorem}

In order to prove a new characterization of cotype in terms of absolutely
summing polynomials we need the following Lemma:

\begin{lemma}
\label{todoponto} If $\mathcal{P}_{as(r;s)(E)}(^{n}E;F)$ $=\mathcal{P}%
(^{n}E;F)$ then $\mathcal{L(}E;F)$ $=\mathcal{L}_{as(r;s)}(E;F).$
\end{lemma}

Proof. (Inspired on the proof of Dvoretzky Rogers Theorem for polynomials
\cite{Matos1})

It is clear that $r\geq s.$ Let us consider a continuous linear mapping $%
T:E\rightarrow F.$ Define an $n$-homogeneous polynomial
\begin{equation*}
P(x)=\varphi(x)^{n-1}T(x)
\end{equation*}
where $\varphi$ is a non null continuous linear functional. Then, choosing $%
a\notin Ker(\varphi),$ we have
\begin{equation*}
dP(a)(x)=(n-1)\varphi(a)^{n-2}\varphi(x)T(a)+\varphi(a)^{n-1}T(x).
\end{equation*}
It is not hard to see that $dP(a)$ is absolutely $(r;s)$-summing (see Matos
\cite{Matos1}) and since $\varphi$ is absolutely $(r;s)$-summing, it follows
that $T$ is absolutely $(r;s)$-summing.$\blacksquare$

It is worth remarking that the converse of Lemma \ref{todoponto} does not
hold. In fact,
\begin{equation*}
\mathcal{L}(l_{2};\mathbb{K)=}\mathcal{L}_{as(2;2)}(l_{2};\mathbb{K)}\text{
and }\mathcal{P}(^{2}l_{2};\mathbb{K)\neq}\mathcal{P}_{as(2;2)}(l_{2};%
\mathbb{K)}.
\end{equation*}
Now we have another characterization of cotype:

\begin{theorem}
Let $n\geq1.$ $E$ has cotype $q>2$ if, and only if,
\begin{equation*}
\mathcal{P}(^{n}E;E)=\mathcal{P}_{as(q;1)(E)}(^{n}E;E).
\end{equation*}
\end{theorem}

Proof. If $\mathcal{P}(^{n}E;E)=\mathcal{P}_{as(q;1)(E)}(^{n}E;E)$ then, by
Lemma \ref{todoponto}, $id:E\rightarrow E$ is $(q;1)$-summing and
consequently $E$ has cotype $q$. Theorem \ref{bbv} yields the converse.$%
\blacksquare$

The following recent Theorem of D.Perez \cite{David}, that generalizes a $2$%
- linear result of Floret-Botelho \cite{Botelho} and Tonge-Melendez\cite
{Tonge}, is an important instrument for other multilinear and holomorphic
results, as we will see later.

\begin{theorem}
(D.Perez \cite{David})\label{Perezz} If each $X_{j}$ is an $\mathcal{L}%
_{\infty,\lambda_{j}}$ space, then every continuous $n$-linear mapping $%
(n\geq2)$ from $X_{1}\times...\times X_{n}$ into $\mathbb{K}$ is $(1;2,...,2)
$-summing at the origin and
\begin{equation*}
\Vert T\Vert_{as(1;2,...,2)}\leq K_{G}3^{\frac{n-2}{2}}\Vert T\Vert
\prod\limits_{j=1}^{n}\lambda_{j}.
\end{equation*}
\end{theorem}

The polynomial version of this Theorem is immediate and will be useful for
us in the last section of this paper.

\begin{corollary}
\label{Perez2}If $X$ is an $\mathcal{L}_{\infty,\lambda}$ space then every
continuous scalar valued $n$-homogeneous polynomial $(n\geq2)$ $%
P:X\rightarrow $ $\mathbb{K}$ is $(1;2)$-summing at the origin and
\begin{equation*}
\Vert P\Vert_{as(1;2)}\leq K_{G}3^{\frac{n-2}{2}}\Vert P\Vert\lambda^{n}.
\end{equation*}
\end{corollary}

We can explore last Theorem and the cotype of the range as follows:

\begin{theorem}
\label{77}If each $X_{j}$ is an $\mathcal{L}_{\infty,\lambda_{j}}$ space and
$F$ has cotype $q\neq\infty,$ then every continuous $n$-linear mapping from $%
X_{1}\times...\times X_{n}$ $\ $into $F$ is $(q;2,...,2)$-summing at the
origin and
\begin{equation*}
\Vert T\Vert_{as(q;2,...,2)}\leq C_{q}(F)K_{G}3^{\frac{n-2}{2}}\Vert
T\Vert\prod\limits_{j=1}^{n}\lambda_{j}.
\end{equation*}
In particular, if $X$ is an $\mathcal{L}_{\infty,\lambda}$ space and $F$ has
cotype $q\neq\infty$, then
\begin{equation*}
\mathcal{P}(^{n}X;F)=\mathcal{P}_{as(q;2)}(^{n}X;F)
\end{equation*}
and
\begin{equation}
\Vert P\Vert_{as(q;2)}\leq C_{q}(F)K_{G}3^{\frac{n-2}{2}}\Vert P\Vert
\lambda^{n}.   \label{referenciar}
\end{equation}
\end{theorem}

Proof. \ Let $(f_{j}^{(1)})_{j=1}^{\infty}\in
l_{2}^{w}(X_{1}),...,(f_{j}^{(n)})_{j=1}^{\infty}\in l_{2}^{w}(X_{n}).$
Since
\begin{equation*}
\mathcal{L}(X_{1},...,X_{n};\mathbb{K})=\mathcal{L}%
_{as(1;2,...,2)}(X_{1},...,X_{n};\mathbb{K})
\end{equation*}
and for every $R\in\mathcal{L}(X_{1},...,X_{n};\mathbb{K})$ we have
\begin{equation*}
\Vert R\Vert_{as(1;2,...,2)}\leq K_{G}3^{\frac{n-2}{2}}\Vert R\Vert
\prod\limits_{j=1}^{n}\lambda_{j}
\end{equation*}
then
\begin{equation*}
(\sum_{j=1}^{\infty}\Vert T(f_{j}^{(1)},...,f_{j}^{(n)})\Vert^{q})^{\frac
{1}{q}}\leq
C_{q}(F)\Vert(T(f_{j}^{(1)},...,f_{j}^{(n)}))_{j=1}^{\infty}\Vert_{w,1}=
\end{equation*}
\begin{equation*}
=C_{q}(F)\sup_{y^{,}\in B_{F^{,}}}\sum_{j=1}^{\infty}\mid(y^{\prime}\circ
T)(f_{j}^{(1)},...,f_{j}^{(n)})\mid=
\end{equation*}
\begin{equation*}
\leq C_{q}(F)\sup_{y^{,}\in B_{F^{,}}}\Vert y^{\prime}\circ T\Vert
_{as(1;2,...,2)}\Vert(f_{j}^{(1)})_{j=1}^{\infty}\Vert_{w,2}...\Vert
(f_{j}^{(n)})_{j=1}^{\infty}\Vert_{w,2}\leq
\end{equation*}
\begin{equation*}
\leq C_{q}(F)\sup_{y^{,}\in B_{F^{,}}}C\Vert y^{\prime}\circ T\Vert\Vert
(f_{j}^{(1)})_{j=1}^{\infty}\Vert_{w,2}...\Vert(f_{j}^{(n)})_{j=1}^{\infty
}\Vert_{w,2}\leq
\end{equation*}
\begin{equation*}
\leq C_{q}(F)\sup_{y^{,}\in B_{F^{,}}}C\Vert y^{\prime}\Vert\Vert T\Vert
\Vert(f_{j}^{(1)})_{j=1}^{\infty}\Vert_{w,2}...\Vert(f_{j}^{(n)})_{j=1}^{%
\infty}\Vert_{w,2}\leq
\end{equation*}
\begin{equation*}
\leq C_{q}(F)C\Vert
T\Vert\Vert(f_{j}^{(1)})_{j=1}^{\infty}\Vert_{w,2}...%
\Vert(f_{j}^{(n)})_{j=1}^{\infty}\Vert_{w,2}
\end{equation*}
where $C=K_{G}3^{\frac{n-2}{2}}\prod\limits_{j=1}^{n}\lambda_{j}.\blacksquare
$

As a consequence of the last Theorem, we obtain a generalization of a
bilinear result of Botelho (\cite{Botelho}), answering a question posed in
\cite{Botelho1}:

\begin{theorem}
\label{102} If $n\geq2$ and each $X_{j}$ is an $\mathcal{L}_{\infty
,\lambda_{j}}$ space then every continuous $n$-linear mapping $%
T:X_{1}\times...\times X_{n}\rightarrow\mathbb{K}$ is $(2;2,...,2,\infty)$%
-summing at the origin and
\begin{equation*}
\Vert T\Vert_{as(2;2,...,2,\infty)}\leq C_{2}(X_{n}%
%TCIMACRO{\UNICODE{0xb4}}%
%BeginExpansion
{\acute{}}%
%EndExpansion
)K_{G}3^{\frac{n-3}{2}}\Vert T\Vert\prod\limits_{j=1}^{n}\lambda_{j}\text{ }%
(\forall n\geq3).
\end{equation*}
\end{theorem}

Proof. \ Let $T:X_{1}\times...\times X_{n}\rightarrow\mathbb{K}$ be a
continuous $n$-linear mapping. Then

$T_{1}:X_{1}\times...\times X_{n-1}\rightarrow X_{n}{}^{\prime}$ is $%
(2;2,...,2)$-summing since $X_{n}{}^{\prime}$ has cotype $2.$ So,
\begin{equation*}
(\sum\limits_{j=1}^{\infty}\Vert T_{1}(x_{j}^{(1)},...,x_{j}^{(n-1)})\Vert
^{2})^{1/2}\leq
C\parallel(x_{j}^{(1)})_{j=1}^{\infty}\parallel_{w,2}...%
\parallel(x_{j}^{(n-1)})_{j=1}^{\infty}\parallel_{w,2}
\end{equation*}
and
\begin{equation*}
(\sum\limits_{j=1}^{\infty}\underset{x_{j}^{(n)}\in B_{X_{n}}}{\sup}\Vert
T_{1}(x_{j}^{(1)},...,x_{j}^{(n-1)})(x_{j}^{(n)})\Vert^{2})^{1/2}\leq
\end{equation*}
\begin{equation*}
\leq C\parallel(x_{j}^{(1)})_{j=1}^{\infty}\parallel_{w,2}...\parallel
(x_{j}^{(n-1)})_{j=1}^{\infty}\parallel_{w,2}.
\end{equation*}
If $(x_{j}^{(n)})_{j=1}^{\infty}\in l_{\infty}(X_{n})$ does not vanish, we
have
\begin{equation*}
(\sum\limits_{j=1}^{\infty}\Vert T_{1}(x_{j}^{(1)},...,x_{j}^{(n-1)})(\frac{%
x_{j}^{(n)}}{\parallel(x_{j}^{(n)})_{j=1}^{\infty}\parallel_{\infty}}%
)\Vert^{2})^{1/2}\leq
\end{equation*}
\begin{equation*}
\leq C\parallel(x_{j}^{(1)})_{j=1}^{\infty}\parallel_{w,2}...\parallel
(x_{j}^{(n-1)})_{j=1}^{\infty}\parallel_{w,2}.
\end{equation*}
Hence
\begin{equation*}
(\sum\limits_{j=1}^{\infty}\Vert
T_{1}(x_{j}^{(1)},...,x_{j}^{(n-1)})(x_{j}^{(n)})\Vert^{2})^{1/2}\leq
\end{equation*}
\begin{equation*}
\leq C\parallel(x_{j}^{(1)})_{j=1}^{\infty}\parallel_{w,2}...\parallel
(x_{j}^{(n-1)})_{j=1}^{\infty}\parallel_{w,2}\parallel(x_{j}^{(n)})_{j=1}^{%
\infty}\parallel_{\infty}
\end{equation*}
and
\begin{equation*}
(\sum\limits_{j=1}^{\infty}\Vert T(x_{j}^{(1)},...,x_{j}^{(n)})\Vert
^{2})^{1/2}\leq
C\parallel(x_{j}^{(1)})_{j=1}^{\infty}\parallel_{w,2}...%
\parallel(x_{j}^{(n-1)})_{j=1}^{\infty}\parallel_{w,2}%
\parallel(x_{j}^{(n)})_{j=1}^{\infty}\parallel_{\infty}
\end{equation*}
where $C=C_{2}(X_{n}%
%TCIMACRO{\UNICODE{0xb4}}%
%BeginExpansion
{\acute{}}%
%EndExpansion
)K_{G}3^{\frac{n-3}{2}}\Vert T\Vert\prod\limits_{j=1}^{n}\lambda_{j}.$

The case $(x_{j}^{(n)})_{j=1}^{\infty}=0$ does not offer any trouble.$%
\blacksquare$

For $n=2$, Theorem \ref{102} has the following version:

\begin{proposition}
\label{apple}If $X$ is an $\mathcal{L}_{\infty}$ space, then every
continuous $2$ linear mapping $T:X\times E\rightarrow\mathbb{K}$ \ with $\cot
$ $E^{\prime}=q=2$ \ is $(r;r,\infty)$-summing at the origin for every $%
r\geq2$. If $\cot$ $E^{\prime}=q>2,$ then $T$ is $(r;r,\infty)$ and $%
(q;p,\infty )$-summing at the origin for every $r>q$ and $p<q.$
\end{proposition}

Proof. (Case $q=2$) Let $T:X\times E\rightarrow\mathbb{K}$ be a continuous
bilinear mapping. Then $T_{1}:X\rightarrow E^{\prime}$ is $(r;r)$-summing
since $E^{\prime}$ has cotype $2$ \cite{Dubinsky}. Hence
\begin{equation*}
(\sum\limits_{j=1}^{\infty}\Vert T_{1}(x_{j})\Vert^{r})^{1/r}\leq
C\parallel(x_{j})_{j=1}^{\infty}\parallel_{w,r}
\end{equation*}
and thus
\begin{equation*}
(\sum\limits_{j=1}^{\infty}\underset{y_{j}\in B_{E}}{\sup}\Vert
T_{1}(x_{j})(y_{j})\Vert^{r})^{1/r}\leq
C\parallel(x_{j})_{j=1}^{\infty}\parallel_{w,r}.
\end{equation*}
If $(y_{j})_{j=1}^{\infty}\in l_{\infty}(E)$ does not vanish, we have
\begin{equation*}
(\sum\limits_{j=1}^{\infty}\Vert T_{1}(x_{j})(\frac{y_{j}}{\parallel
(y_{j})_{j=1}^{\infty}\parallel_{\infty}})\Vert^{r})^{1/r}\leq C\parallel
(x_{j})_{j=1}^{\infty}\parallel_{w,r}.
\end{equation*}
Hence
\begin{equation*}
(\sum\limits_{j=1}^{\infty}\Vert T_{1}(x_{j})(y_{j})\Vert^{r})^{1/r}\leq
C\parallel(x_{j})_{j=1}^{\infty}\parallel_{w,r}\parallel(y_{j})_{j=1}^{%
\infty }\parallel_{\infty}
\end{equation*}
and
\begin{equation*}
(\sum\limits_{j=1}^{\infty}\Vert T(x_{j},y_{j})\Vert^{r})^{1/r}\leq
C\parallel(x_{j})_{j=1}^{\infty}\parallel_{w,r}\parallel(y_{j})_{j=1}^{%
\infty }\parallel_{\infty}.
\end{equation*}
The case$\ (y_{j})_{j=1}^{\infty}=0$ does not offer any problem.

A linear result of Maurey (see \cite{Diestel}, page 223, Th. 11.14 a)
provides, through the same reasoning, the proof of the case $q>2$.$%
\blacksquare$

Applying the same ideas we have the statement below:

\begin{proposition}
\label{20}If each $X_{j}$ is an $\mathcal{L}_{\infty}$ space$,$ then every
continuous $n$-linear mapping $T:X_{1}\times...\times X_{n}\times
E\rightarrow\mathbb{K}$ \ with $\cot$ $E^{\prime}=q\geq2,q\neq\infty$ is $%
(q;2,...,2,\infty)$-summing at the origin.
\end{proposition}

Theorem \ref{102} can also be used to obtain other results. For example:

\begin{theorem}
\label{00}If each $X_{j}$ is an $\mathcal{L}_{\infty}$ space and $%
T:X_{1}\times...\times X_{n}\rightarrow\mathbb{K}$ $\ $is a continuous $n$%
-linear mapping, then

$n=2\Rightarrow T$ is $(r;r,r)$-summing on $X_{1}\times X_{2},$ for every $%
r\geq2$.

$n\geq3\Rightarrow T$ is $(r;2,...,2,r)$-summing on $X_{1}\times...\times
X_{n}$ \ for every $r\geq2.$
\end{theorem}

Proof. The case $n=2$ is the easiest and we will omit the proof. For the
case $n=3,$ let $(x_{j})_{j=1}^{\infty}\in
l_{2}^{w}(X_{1}),(y_{j})_{j=1}^{\infty }\in l_{2}^{w}(X_{2})$ and $%
(z_{j})_{j=1}^{\infty}\in l_{2}^{w}(X_{3}).$ Then
\begin{equation*}
(\sum\limits_{j=1}^{\infty}\Vert T(a+x_{j},b+y_{j},c+z_{j})-T(a,b,c)\Vert
^{r})^{\frac{1}{r}}=
\end{equation*}
\begin{equation*}
=(\sum\limits_{j=1}^{\infty}(\Vert T(a,y_{j},z_{j})\Vert^{r})^{\frac{1}{r}%
}+(\sum\limits_{j=1}^{\infty}(\Vert T(x_{j},b,c)\Vert^{r})^{\frac{1}{r}}+
\end{equation*}
\begin{equation*}
+(\sum\limits_{j=1}^{\infty}(\Vert T(x_{j},y_{j},c)\Vert^{r})^{\frac{1}{r}%
}+(\sum\limits_{j=1}^{\infty}(\Vert T(x_{j},b,z_{j})\Vert^{r})^{\frac{1}{r}%
}+(\sum\limits_{j=1}^{\infty}(\Vert T(a,b,z_{j})\Vert^{r})^{\frac{1}{r}}
\end{equation*}
\begin{equation*}
+(\sum\limits_{j=1}^{\infty}(\Vert T(a,y_{j},c)\Vert^{r})^{\frac{1}{r}%
}+(\sum\limits_{j=1}^{\infty}(\Vert T(x_{j},y_{j},z_{j})\Vert^{r})^{\frac{1}{%
r}}<\infty
\end{equation*}
since every linear mapping is $(r;r)$ and $(r;2)$-summing, every such
bilinear mapping is $(r;r,r)$, $(r;2,r)$ and $(r;2,2)$-summing and every
such $3$-linear mapping above is $(r;2,2,r)$-summing at the origin.

For $n>3$ we use an inductive principle.

Theorem \ref{77} can be extended as follows:

\begin{theorem}
\label{200} If each $X_{j}$ is an $\mathcal{L}_{\infty}$ space and $\cot F=q,
$ then every continuous $n$-linear mapping from $X_{1}\times...\times X_{n}$
into $F$ is $(q;2,...,2)$-summing on $X_{1}\times...\times X_{n}.$
\end{theorem}

Proof. If $q=2,$ it is enough to use the last reasoning with Theorem \ref{77}
and the Dubinsky-Pe\l czy\'{n}ski-Rosenthal (\cite{Diestel} page 223, Th.
11.14 (a) or \cite{Dubinsky}) result which asserts that every linear mapping
from an $\mathcal{L}_{\infty}$ space into $F$( with $\cot F=2$ ) is \ $(2;2)$%
-summing.

If $q>2,$ we shall use the same reasoning with the Maurey (\cite{Diestel}
page 223, Th. 11.14(b)) result which asserts that every linear mapping from
an $\mathcal{L}_{\infty}$ space into $F$( with $\cot$ $F=q>2$ ) is $(q;p)$%
-summing for each $p<q.\blacksquare$

\section{r-fully absolutely summing multilinear mappings}

The following definition is inspired in the work of Matos \cite{fully} which
is being developed by M.L.V. Souza in his doctoral dissertation.

\begin{definition}
A continuous $n$-linear mapping $T:E_{1}\times...\times E_{n}\rightarrow F$
will be said $r$-fully $(p;q_{1},...,q_{n})$-summing if
\begin{equation*}
\sum_{j_{1},...,j_{r}=1}^{\infty}\left\|
T(x_{j_{1}}^{(1)},...,x_{j_{r}}^{(r)},x_{j_{r}}^{(r+1)},...,x_{j_{r}}^{(n)}%
\right\| ^{p}<\infty
\end{equation*}
whenever $(x_{k}^{(l)})_{k=1}^{\infty}\in l_{q_{l}}^{w}(E_{l}),l=1,...,n$.
In this case we will write
\begin{equation*}
T\in\mathcal{L}_{f(r)as(p;q_{1},...,q_{n})}(E_{1},...,E_{n};F).
\end{equation*}
\end{definition}

When $r=1$, we have the $(p;q_{1},...,q_{n})$-summing mappings and when $r=n$
we call $T$ just by fully $(p;q_{1},...,q_{n})$-summing, which is a concept
introduced by Matos \cite{fully}.

A natural question is: Does every $n$-linear mappings from $\mathcal{L}%
_{\infty}$ spaces into $\mathbb{K}$ is fully $(2;2,...,2)$-summing?

We will show in Corollary \ref{tres} that Theorem \ref{102} give us partial
answers.

\begin{theorem}
If $\mathcal{L}(E_{n};F)=\mathcal{L}_{as(q;r)}(E_{n};F),$ then
\begin{equation*}
\mathcal{L}_{as(q,p_{1},...,p_{n-1},\infty)}(E_{1},...,E_{n};F)\subset
\mathcal{L}_{f(2)as(q;p_{1},...,p_{n-1},r)}(E_{1},...,E_{n};F).
\end{equation*}
\end{theorem}

Proof. Let us consider $T\in\mathcal{L}_{as(q;p_{1},...,p_{n-1},%
\infty)}(E_{1},...,E_{n};F).$ If
\begin{equation*}
(x_{k}^{(1)})_{k=1}^{\infty}\in
l_{p_{1}}^{w}(E_{1}),...,(x_{k}^{(n-1)})_{k=1}^{\infty}\in
l_{p_{n-1}}^{w}(E_{n-1}),(y_{k})_{k=1}^{\infty}\in l_{r}^{w}(E),
\end{equation*}
then, for each $k$ fixed,
\begin{equation*}
\sum\limits_{j=1}^{\infty}\left\|
T(x_{k}^{(1)},...,x_{k}^{(n-1)},y_{j})\right\| ^{q}\leq\left\|
(y_{j})_{j=1}^{\infty}\right\| _{w,r}^{q}\left\|
T(x_{k}^{(1)},...,x_{k}^{(n-1)},.)\right\| _{as(q,r)}^{q}\leq
\end{equation*}
\begin{equation*}
\leq\left\| (y_{j})_{j=1}^{\infty}\right\| _{w,r}^{q}C\left\|
T(x_{k}^{(1)},...,x_{k}^{(n-1)},.)\right\| ^{q}\leq
\end{equation*}
\begin{equation*}
\leq\left\| (y_{j})_{j=1}^{\infty}\right\| _{w,r}^{q}C(\left\|
T(x_{k}^{(1)},...,x_{k}^{(n-1)},z_{k})\right\| ^{q}+\frac{1}{2^{k}}).
\end{equation*}
where each $z_{k}$ belongs to the unit ball $B_{E_{n}}.$

Therefore

\begin{equation*}
\sum\limits_{k=1}^{\infty}\sum\limits_{j=1}^{\infty}\left\|
T(x_{k}^{(1)},...,x_{k}^{(n-1)},y_{j})\right\| ^{q}\leq
\end{equation*}
\begin{equation*}
\leq\left\| (y_{j})\right\| _{w,r}^{q}C\sum\limits_{k=1}^{\infty}(\left\|
T(x_{k}^{(1)},...,x_{k}^{(n-1)},z_{k})\right\| ^{p}+\frac{1}{2^{k}}%
)<\infty.\blacksquare
\end{equation*}

\begin{corollary}
\label{tres}If each $E_{k}$ is an $\mathcal{L}_{\infty}$ space, we have
\begin{equation*}
\mathcal{L}(E_{1},...,E_{n};\mathbb{K)}=\mathcal{L}%
_{f(2)as(2;2,...,2,2)}(E_{1},...,E_{n};\mathbb{K}).
\end{equation*}
\end{corollary}

Proof. It suffices to realize that $\mathcal{L}(E_{n};\mathbb{K})=\mathcal{L}%
_{as(2;2)}(E_{n};\mathbb{K})$ $\ $and apply last Theorem and Theorem \ref
{102}.$\blacksquare$

\section{Other results}

An important and broadly used result is the Generalized H\"{o}lder%
%TCIMACRO{\UNICODE{0xb4}}%
%BeginExpansion
\'{}%
%EndExpansion
s Inequality, which is a natural instrument to deal with absolutely summing
multilinear mappings.

\begin{theorem}[Generalized H\"{o}lder%
%TCIMACRO{\UNICODE{0xb4}}%
%BeginExpansion
\'{}%
%EndExpansion
s Inequality]
If $\frac{1}{p}\leq\frac{1}{p_{1}}+...+\frac{1}{p_{n}},$ then
\begin{equation*}
(\sum_{j=1}^{\infty}\mid a_{j}^{(1)}...a_{j}^{(n)}\mid^{p})^{\frac{1}{p}%
}\leq(\sum_{j=1}^{\infty}\mid a_{j}^{(1)}\mid^{p_{1}})^{\frac{1}{p_{1}}%
}...(\sum_{j=1}^{\infty}\mid a_{j}^{(n)}\mid^{p_{n}})^{\frac{1}{p_{n}}}.
\end{equation*}
\end{theorem}

If $T:E_{1}\times...\times E_{n}\rightarrow F$ is a continuous multilinear
mapping where at least one of the spaces which compose the Banach spaces of
the domain has finite cotype, we can state the following result.

\begin{theorem}
\label{cotipodom}If $T:E_{1}\times...\times E_{n}\rightarrow F$ is a
continuous multilinear mapping, $q_{j}=\cot E_{j},$ $j=1,...,n,$ and at
least one of the $q_{j}$ finite, then, for any choice of $a_{j}\in\lbrack
q_{j},\infty],$ with at least one of the $a_{j}$ finite, $T$ is $%
(s;b_{1},...,b_{n})$-summing at the origin, for any $s>0,$ such that $\frac{1%
}{s}\leq\frac{1}{a_{1}}+...+\frac{1}{a_{n}},$ with $b_{j}=1,$ if $%
a_{j}<\infty,$ and $b_{j}=\infty$ if $a_{j}=\infty.$
\end{theorem}

Proof. Obvious, using Theorem \ref{Shera}, after some reasoning on how to
optimize the use of the Generalized H\"{o}lder%
%TCIMACRO{\UNICODE{0xb4}}%
%BeginExpansion
\'{}%
%EndExpansion
s Inequality.

As a corollary, we have the a result due to Botelho \cite{Botelho}.

\begin{corollary}
\label{Bot} If $T:E_{1}\times...\times E_{n}\rightarrow F$ is a continuous
multilinear mapping and $q_{j}=\cot E_{j}<\infty$ for every $j=1,...,n$,
then $T$ is $(s;1,...,1)$-summing at the origin for any $s>0$ such that $%
\frac {1}{s}\leq\frac{1}{q_{1}}+...+\frac{1}{q_{n}}.$
\end{corollary}

Theorem \ref{cotipodom} shows that even if just one of the spaces of the
domain has finite cotype, the multilinear mapping is still well behaved. As
an illustration we can see the example below.

\begin{example}
If $E$ has cotype $p$, then every $T:C(K)\times...\times C(K)\times
E\rightarrow F$ is $(p;\infty,...,\infty,1)$-summing at the origin.
\end{example}

The following results show more about the mechanism of absolutely summing
mappings.

\begin{proposition}
\label{mh}If $\mathcal{L}(E_{1},...,E_{n};F)=\mathcal{L}%
_{as(r;s_{1},...,s_{t},\infty,...,\infty)}(E_{1},...,E_{n};F)$ then
\begin{equation*}
\mathcal{L}(E_{1},...,E_{t};F)=\mathcal{L}%
_{as(r;s_{1},...,s_{t})}(E_{1},...,E_{t};F).
\end{equation*}
\end{proposition}

Proof. Given $T$ $\in$ $\mathcal{L}(E_{1},...,E_{t};F)$ let us define
\begin{equation*}
S(a_{1},...,a_{n})=T(a_{1},...,a_{t})\varphi_{t+1}(a_{t+1})...\varphi
_{n}(a_{n})
\end{equation*}
where $\varphi_{t+1},...,\varphi_{n}$ are non trivial bounded linear
functionals. Let $b_{t+1},...,b_{n}$ be such that
\begin{equation*}
\varphi_{t+1}\left( b_{t+1}\right) =...=\varphi_{n}(b_{n})=1.
\end{equation*}
It follows that $T\in\mathcal{L}_{as(r;s_{1},...,s_{t})}(E_{1},...,E_{t};F).$
In fact, if $(x_{j}^{l})_{j=1}^{\infty}\in l_{s_{l}}^{w}(E_{l})$ we have
\begin{equation*}
\sum_{j=1}^{\infty}\Vert T(x_{j}^{(1)},...,x_{j}^{(t)})\Vert^{r}=\sum
_{j=1}^{\infty}\Vert S(x_{j}^{(1)},...,x_{j}^{(t)},b_{t+1},...,b_{n})\Vert
^{r}<\infty.\blacksquare
\end{equation*}
The next statement suggested by Matos extend the Lemma 3.2 of \cite{Botelho}:

\begin{proposition}
\label{Lemma} If $\mathcal{L}(E_{1},...,E_{n};F)=\mathcal{L}%
_{as(r;s_{1},...,s_{t},\infty,...,\infty)}(E_{1},...,E_{n};F),$ then
\begin{equation*}
\mathcal{L}(E_{1},...,E_{t};\mathcal{L}(E_{t+1},...,E_{n};F))=\mathcal{L}%
_{as(r;s_{1},...,s_{t})}(E_{1},...,E_{t};\mathcal{L}(E_{t+1},...,E_{n};F))
\end{equation*}
and conversely.
\end{proposition}

Proof: Suppose $\ $%
\begin{equation*}
\mathcal{L}(E_{1},...,E_{n};F)=\mathcal{L}_{as(r;s_{1},...,s_{t},\infty,...,%
\infty)}(E_{1},...,E_{n};F).
\end{equation*}
Let $T:E_{1}\times...\times E_{t}\longrightarrow\mathcal{L}%
(E_{t+1},...,E_{n};F)$ be a continuous multilinear mapping. We have
\begin{equation*}
(\sum_{j=1}^{\infty}\Vert T(x_{1}^{(j)},...,x_{t}^{(j)})\Vert^{r})^{\frac
{1}{r}}=(\sum_{j=1}^{\infty}\underset{\underset{k=t+1,...,n}{\Vert
y_{k}\Vert\leq1}}{\sup}\Vert
T(x_{1}^{(j)},...,x_{t}^{(j)})(y_{t+1},...,y_{n})\Vert^{r})^{\frac{1}{r}%
}\leq
\end{equation*}
\begin{equation*}
\leq(\sum_{j=1}^{\infty}\Vert
T(x_{1}^{(j)},...,x_{t}^{(j)})(y_{t+1}^{(j)},...,y_{n}^{(j)})\Vert^{r}+\frac{%
1}{2^{j}})^{\frac{1}{r}}=
\end{equation*}
\begin{equation*}
=(\sum_{j=1}^{\infty}\Vert
T_{0}(x_{1}^{(j)},...,x_{t}^{(j)},y_{t+1}^{(j)},...,y_{n}^{(j)})\Vert^{r}+%
\frac{1}{2^{j}})^{\frac{1}{r}}<\infty
\end{equation*}
if $(x_{1}^{(j)})\in l_{s_{1}}^{w}(E_{1}),...,(x_{t}^{(j)})\in
l_{s_{t}}^{w}(E_{t}).$

On the other hand, suppose
\begin{equation*}
\mathcal{L}(E_{1},...,E_{t};\mathcal{L}(E_{t+1},...,E_{n};F)=\mathcal{L}%
_{as(r;s_{1},...,s_{t})}(E_{1},...,E_{t};\mathcal{L}(E_{t+1},...,E_{n};F).
\end{equation*}
If $T:E_{1}\times...\times E_{n}\longrightarrow F,$ and
\begin{equation*}
(x_{1}^{(j)})_{j=1}^{\infty}\in
l_{s_{1}}^{w}(E_{1}),...,(x_{t}^{(j)})_{j=1}^{\infty}\in
l_{s_{t}}^{w}(E_{t}),(y_{t+1}^{(j)})_{j=1}^{\infty}\in
l_{\infty}(E_{t}),...,(y_{n}^{(j)})_{j=1}^{\infty}\in l_{\infty}(E_{n}),
\end{equation*}
we have
\begin{equation*}
\sum_{j=1}^{\infty}\Vert
T(x_{1}^{(j)},...,x_{t}^{(j)},y_{t+1}^{(j)},...,y_{n}^{(j)})\Vert^{r})^{%
\frac{1}{r}}=(\sum_{j=1}^{\infty}\Vert
T_{1}(x_{1}^{(j)},...,x_{t}^{(j)})(y_{t+1}^{(j)},...,y_{n}^{(j)})\Vert
^{r})^{\frac{1}{r}}\leq
\end{equation*}
\begin{equation}
\leq\Vert(y_{t+1}^{(j)})\Vert_{\infty}...\Vert(y_{n}^{(j)})\Vert_{\infty}%
\sum_{j=1}^{\infty}\Vert T_{1}(x_{1}^{(j)},...,x_{t}^{(j)})\Vert^{r})^{\frac{%
1}{r}}<\infty.   \label{dssa}
\end{equation}
We can see that it is also true that
\begin{equation*}
T\in\mathcal{L}_{as(r;s_{1},...,s_{t},\infty,...,\infty)}(E_{1},...,E_{n};F)%
\Rightarrow
\end{equation*}
\begin{equation*}
\Rightarrow T_{1}\in\mathcal{L}_{as(r;s_{1},...,s_{t},)}(E_{1},...,E_{t};%
\mathcal{L}(E_{t+1},...,E_{n};F))
\end{equation*}
and
\begin{equation*}
T\in\mathcal{L}_{as(r;s_{1},...,s_{t},)}(E_{1},...,E_{t};\mathcal{L}%
(E_{t+1},...,E_{n};F))\Rightarrow
\end{equation*}
\begin{equation*}
\Rightarrow T_{0}\in\mathcal{L}_{as(r;s_{1},...,s_{t},\infty,...,\infty
)}(E_{1},...,E_{n};F).\blacksquare
\end{equation*}

\begin{remark}
The reader shall note that the converse of Proposition \ref{mh} cannot hold.
In fact, we know that $\mathcal{L}(E;\mathbb{K})=\mathcal{L}_{as(1;1)}(E;%
\mathbb{K}).$ If \ the converse of Proposition \ref{mh} held, we would have $%
\mathcal{L}(E,E;\mathbb{K})=\mathcal{L}_{as(1;1,\infty)}(E,E;\mathbb{K})$
and by Proposition \ref{Lemma}
\begin{equation*}
\mathcal{L}(E;E%
%TCIMACRO{\UNICODE{0xb4}}%
%BeginExpansion
{\acute{}}%
%EndExpansion
)=\mathcal{L}_{as(1;1)}(E;E%
%TCIMACRO{\UNICODE{0xb4}}%
%BeginExpansion
{\acute{}}%
%EndExpansion
)
\end{equation*}
which is impossible, in general (see \cite{Lindenstrauss}).
\end{remark}

Proposition \ref{Lemma} also furnishes an Inclusion Theorem for bilinear
mappings.

\begin{proposition}
\label{telnet}( Inclusion for bilinear mappings)

If $r>s$ then $\mathcal{L}_{as(s;s,\infty)}(E_{1},E_{2};F)\subset \mathcal{L}%
_{as(r;r,\infty)}(E_{1},E_{2};F).$
\end{proposition}

Proof. If $r>s$ and $T\in\mathcal{L}_{as(s;s,\infty)}(E_{1},E_{2};F),$ then
by Proposition \ref{Lemma}, $T_{1}:E_{1}\rightarrow\mathcal{L}(E_{2};F)$ is $%
(s;s)$-summing. By the Inclusion Theorem for linear mappings, $T_{1}$ will
be $(r;r)$-summing and again by the Proposition \ref{Lemma}, $T$ will be $%
(r;r,\infty)$-summing at the origin.$\blacksquare$

\begin{example}
The famous Grothendieck%
%TCIMACRO{\UNICODE{0xb4}}%
%BeginExpansion
\'{}%
%EndExpansion
s Theorem, which asserts that every linear operator from an $\mathcal{L}_{1}$
space into an $\mathcal{L}_{2}$ space is $(1;1)$-summing, and Proposition
\ref{Lemma} lead us to conclude that if $E_{1}$ and $E_{2}$ are $\mathcal{L}%
_{1}$ and $\mathcal{L}_{2}$ spaces respectively, then
\begin{equation*}
\mathcal{L}(E_{1},E_{2};\mathbb{K})=\mathcal{L}_{as(1;1,\infty)}(E_{1},E_{2};%
\mathbb{K}).
\end{equation*}
Thus, Proposition \ref{telnet} yields
\begin{equation*}
\mathcal{L}(E_{1},E_{2};\mathbb{K})=\mathcal{L}_{as(r;r,\infty)}(E_{1},E_{2};%
\mathbb{K})
\end{equation*}
for every $r\geq1$. However, despite Grothendieck%
%TCIMACRO{\UNICODE{0xb4}}%
%BeginExpansion
\'{}%
%EndExpansion
s Theorem we know that
\begin{equation*}
\mathcal{L}(l_{1},l_{1};l_{2})\neq\mathcal{L}_{as(1;1,%
\infty)}(l_{1},l_{1};l_{2})
\end{equation*}
and furthermore
\begin{equation*}
\mathcal{L}(l_{1},l_{1};\mathbb{K})\neq\mathcal{L}_{as(1;1,%
\infty)}(l_{1},l_{1};\mathbb{K}).
\end{equation*}
\end{example}

The result below has the same spirit of the last Proposition.

\begin{proposition}
\label{poyt}\ If $T:E_{1}\times...\times E_{n}\rightarrow F$ is $p$%
-dominated, then $T$ is $(\frac{r}{n-1};r,...,r,\infty)$-summing for every $%
r\geq p.$
\end{proposition}

Proof. If $T$ $:E_{1}\times...\times E_{n}\rightarrow F$ is $p$-dominated,
then, by Grothendieck-Pietsch domination Theorem, if $T_{1}:E_{1}\times...%
\times E_{n-1}\rightarrow\mathcal{L}(E_{n};F)$ is such that $T_{1}=\Psi(T),$
we obtain, for $r\geq p,$%
\begin{equation*}
\Vert T_{1}(x_{1},...,x_{n-1})\Vert=\sup_{\Vert y\Vert\leq1}\Vert
T_{1}(x_{1},...,x_{n-1})(y)\Vert=
\end{equation*}
\begin{equation*}
=\sup_{\Vert y\Vert\leq1}\Vert T(x_{1},...,x_{n-1},y)\Vert\leq
\end{equation*}
\begin{equation*}
\leq\sup_{\Vert
y\Vert\leq1}C(\int_{B_{E_{1}^{\prime}}}\mid\varphi(x_{1})\mid^{r}d\mu_{1})^{%
\frac{1}{r}}...(\int_{B_{E_{n}^{\prime}}}\mid \varphi(y)\mid^{r}d\mu_{n})^{%
\frac{1}{r}}\leq
\end{equation*}
\begin{equation*}
\leq C(\int_{B_{E_{1}^{\prime}}}\mid\varphi(x_{1})\mid^{r}d\mu_{1})^{\frac
{1}{r}}...(\int_{B_{E_{n-1}^{\prime}}}\mid\varphi(x_{n-1})\mid^{r}d\mu
_{n-1})^{\frac{1}{r}}.
\end{equation*}
Thus, $T_{1}$ is $r$-dominated and, by Proposition \ref{Lemma}, $%
T=(T_{1})_{0}$ \newline
is $(\frac{r}{n-1}r;r,...,r,\infty)$-summing.

\begin{corollary}
If every $T:E_{1}\times...\times E_{n}\rightarrow F$ is $p$-dominated, then
every $T:E_{j_{1}}\times...\times E_{j_{r}}\rightarrow F,$ with $1\leq r\leq
n$ and $j_{1},...,j_{r}\in\{1,...,n\}$ mutually disjoint, is $p-$dominated.
\end{corollary}

Proof. By Proposition \ref{poyt}, we have
\begin{equation*}
\mathcal{L}(E_{1},...,E_{n};F)=\mathcal{L}_{as(\frac{p}{n-1};p,...,p,\infty
)}(E_{1},...,E_{n};F)
\end{equation*}
and by Proposition \ref{mh} we obtain
\begin{equation*}
\mathcal{L}(E_{1},...,E_{n-1};F)=\mathcal{L}_{as(\frac{p}{n-1}%
;p,...,p)}(E_{1},...,E_{n-1};F).
\end{equation*}
The other cases use the same arguments.$\blacksquare$ \newline
Similar reasoning furnishes the next Corollary.

\begin{corollary}
\label{frew}If every $T:E_{1}\times...\times E_{n}\rightarrow F$ is $p$%
-dominated, then for every permutation $\pi:\{1,...,n\}\rightarrow
\{1,...,n\}$ we have
\begin{equation*}
\mathcal{L}(E_{\pi(1)},...,E_{\pi(t)};\mathcal{L}(E_{\pi(t+1)},...,E_{\pi
(n)};F))=
\end{equation*}
\begin{equation*}
=\mathcal{L}_{as(\frac{p}{t},p,...,p)}(E_{\pi(1)},...,E_{\pi(t)};\mathcal{L}%
(E_{\pi(t+1)},...,E_{\pi(n)};F)).
\end{equation*}
\end{corollary}

The next result is essentially due to Botelho \cite{Botelho}.

\begin{corollary}
If some $E_{j}$ is an $\mathcal{L}_{\infty}$ space, at least one other $E_{k}
$ is infinite dimensional and $\dim F=\infty$, then, regardless of the $%
p\geq1,$ we have
\begin{equation*}
\mathcal{L}(E_{1},...,E_{n};F)\neq\mathcal{L}_{as(\frac{p}{n}%
;p)}(E_{1},...,E_{n};F).
\end{equation*}
\end{corollary}

Proof. There is no loss of generality in assuming $j=1.$ If the equality
held we would have
\begin{equation*}
\mathcal{L}(E_{1};\mathcal{L}(E_{2},...,E_{n};F))=\mathcal{L}%
_{as(p;p)}(E_{1};\mathcal{L}(E_{2},...,E_{n};F))
\end{equation*}
which is a contradiction since $\mathcal{L}(E_{2},...,E_{n};F)$ has only
infinite cotype (see \cite{Botelho},\cite{Dineen}).

\section{Extrapolation Theorems}

The linear theory of absolutely summing operators has some strong
coincidence Theorems (see \cite{Diestel}). Many of them have their
polynomial versions (see \cite{Matos},\cite{Tonge}). We will give a
polynomial and a multilinear version for the Maurey Extrapolation Theorem:

\begin{theorem}[Polynomial Extrapolation Theorem]
\label{ext}If $1<r<p<\infty$ and $X$ is a Banach space such that
\begin{equation}
\mathcal{P}_{as(\frac{p}{n};p)}(^{n}X;l_{p})=\mathcal{P}_{as(\frac{r}{n}%
;r)}(^{n}X;l_{p})   \label{cc}
\end{equation}
then, for every Banach space $Y$ we have
\begin{equation}
\mathcal{P}_{as(\frac{p}{n};p)}(^{n}X;Y)=\mathcal{P}_{as(\frac{1}{n}%
;1)}(^{n}X;Y)\text{ }
\end{equation}
\end{theorem}

Proof. Consider
\begin{equation*}
\varphi:X\rightarrow C(B_{X^{\ast}}):x\mapsto f_{x}
\end{equation*}
where $f_{x}(x^{\ast})=<x^{\ast},x>.$ We will denote $K=B_{X^{\ast}}.$ Let
us denote by $P(K)$ the set of \ all probability measures on $K$ with the
weak star topology. For each $\mu\in P(K)$ define
\begin{equation*}
j_{\mu}:X\subset C(K)\rightarrow L_{p}(\mu)
\end{equation*}
as the restriction of the canonical inclusion from $C(K)$ into $L_{p}(\mu).$

Let $R:X\rightarrow Y$ be an $n$-homogeneous $(\frac{p}{n};p)$-summing
polynomial. The polynomial version of Grothendieck-Pietsch Domination
Theorem tells us that there exists $\mu_{0}\in P(K)$ such that
\begin{equation*}
\Vert Rx\Vert\leq C[\int_{K}\mid<\varphi,x>\mid^{p}d\mu_{0}(\varphi )]^{%
\frac{n}{p}}=C[\int_{K}\mid j_{\mu_{0}}(x)(\varphi)\mid^{p}d\mu
_{0}(\varphi)]^{\frac{n}{p}}=
\end{equation*}
\begin{equation*}
=C\Vert j_{\mu_{0}}(x)\Vert_{L_{p}(\mu_{0})}^{n}\text{ \ \ for every }x\text{
in }X.
\end{equation*}
We must find $\lambda\in P(K)$ and a constant $D$ (depending on $X$ ) such
that
\begin{equation}
\Vert j_{\mu_{0}}(x)\Vert_{L_{p}(\mu_{0})}\leq D\Vert
j_{\lambda}(x)\Vert_{L_{1}(\lambda)}\text{ }\forall x\in X,   \label{ext1}
\end{equation}
and then the Theorem will be proved. Indeed, we will have

\begin{equation*}
\Vert Rx\Vert\leq C\Vert j_{\mu_{0}}(x)\Vert_{L_{p}(\mu_{0})}^{n}\leq
CD\Vert j_{\lambda}(x)\Vert_{L_{1}(\lambda)}^{n}=
\end{equation*}
\begin{equation*}
=C_{1}[\int_{K}\mid j_{\lambda}(x)(x^{\ast})\mid
d\lambda(\varphi)]^{n}=C_{1}[\int_{K}\mid x^{\ast}(x)\mid
d\lambda(\varphi)]^{n}
\end{equation*}
and the Grothendieck-Pietsch Polynomial Domination Theorem yields that $R$
is $(\frac{1}{n};1)$-summing.

In order to prove (\ref{ext1}) it is enough to note that
\begin{equation*}
\mathcal{P}_{as(\frac{p}{n};p)}(^{n}X;l_{p})=\mathcal{P}_{as(\frac{r}{n}%
;r)}(^{n}X;l_{p})
\end{equation*}
imply $\mathcal{L}_{as,p}(X;l_{p})=\mathcal{L}_{as,r}(X;l_{p}),$ and it is
enough to end the proof (this is done in the proof of the linear
Extrapolation Theorem. See Th. 3.17 of \cite{Diestel}).$\blacksquare$

For the multilinear version, the same reasoning give us the following
statement:

\begin{theorem}
If $1<r<p<\infty$ and $X$ is a Banach space such that
\begin{equation*}
\mathcal{L}_{as(\frac{p}{n};p)}(^{n}X;l_{p})=\mathcal{L}_{as(\frac{r}{n}%
;r)}(^{n}X;l_{p})\text{ }
\end{equation*}
then, for every Banach space $Y$, we have
\begin{equation*}
\mathcal{L}_{as(\frac{p}{n};p)}(^{n}X;Y)=\mathcal{L}_{as(\frac{1}{n}%
;1)}(^{n}X;Y).\text{ }
\end{equation*}
\end{theorem}

\section{Absolutely summing mappings}

The concept of absolutely summing mapping (non necessarily multilinear or
polynomial) and the first results and examples are due to M.Matos \cite
{Matos2}.

\begin{definition}
(Matos) A mapping $f:E\rightarrow F$ is absolutely $(s;r)$-summing at $a\in
E $ if $(f(a+x_{j})-f(a))_{j=1}^{\infty}\in l_{s}(F)$ whenever $%
(x_{j})_{j=1}^{\infty}\in l_{r}^{u}(E).$ We say that $f:E\rightarrow F$ is
weakly absolutely $(s,r)$-summing at $a\in E$ if $(f(a+x_{j})-f(a))_{j=1}^{%
\infty}\in l_{s}^{w}(F)$ whenever $(x_{j})_{j=1}^{\infty}\in l_{r}^{u}(E).$
\end{definition}

Since for every $(x_{j})_{j=1}^{\infty}\in l_{r}^{u}(E)$ we have $%
\lim_{m\rightarrow\infty}\left\| (x_{j})_{j=m}^{\infty}\right\| _{w,r}=0,$
it is clear that $\lim_{m\rightarrow\infty}\left\| x_{m}\right\| =0.$
Therefore, there is no loss of generality if, in the definition above, we
restrict ourselves to $(x_{j})_{j=1}^{\infty}\in l_{r}^{u}(E)$ with $\left\|
x_{j}\right\| <\delta$ for all $j$ and some $\delta.$

It is possible to prove that if $f$ $:E\rightarrow F$ is absolutely $(s;r)$%
-summing at $a\in E$ then $f$ is continuous at $a$ \cite{Matos}. The
behavior of $f$ outside an open neighborhood of $a$ is completely irrelevant.

In \cite{Botelho1}, Botelho proves for the complex case, using Cauchy
integral formulas, that, if $\cot E=q,$ every holomorphic entire mapping $%
f:E\rightarrow F$ such that $f(0)=0$ is $(q;1)$-summing at the origin. We
will prove that Cauchy integral formulas are not essential and this result
still holds for the real case and for non zero points.

\begin{lemma}
\label{LLema}If $g:E\rightarrow F$ is analytic at $a\in E$, then there
exists $\delta>0$ such that
\begin{equation*}
\Vert(g(a+x_{j})-g(a))_{j=1}^{\infty}\Vert_{w,1}\leq
D\Vert(x_{j})_{j=1}^{\infty}\Vert_{w,1}\text{ }
\end{equation*}
whenever $\Vert(x_{j})_{j=1}^{\infty}\Vert_{w,p}<\delta.$
\end{lemma}

Proof. If $g:E\rightarrow F$ is analytic at $a$ and $C,c>0$ are such that
\begin{equation*}
\Vert\frac{1}{k!}\overset{\wedge}{d^{k}}g(a)\Vert\leq Cc^{k}\text{ for every
}k
\end{equation*}
then, for each $\varphi\in F`$, we have
\begin{equation*}
\Vert\frac{1}{k!}\overset{\wedge}{d^{k}}\varphi g(a)\Vert=\Vert\varphi\frac
{1}{k!}\overset{\wedge}{d^{k}}g(a)\Vert\leq Cc^{k}\Vert\varphi\Vert\text{
for all }k
\end{equation*}
and hence, by a result of Defant and Voigt (see \cite{Matos}, Theorem 1.6 or
\cite{Matos2}),
\begin{equation*}
\Vert\frac{1}{k!}\overset{\wedge}{d^{k}}\varphi g(a)\Vert_{as(1;1)}\leq
e^{k}Cc^{k}\Vert\varphi\Vert.
\end{equation*}
Let us denote by $\epsilon_{a}>0$ the radius of convergence of $g$ around $%
a. $Thus, if $\Vert(x_{j})_{j=1}^{\infty}\Vert_{w,1}\leq\delta=\min\{\frac
{1}{2ec},\epsilon_{a}\}$\ we can write
\begin{equation*}
\sum_{j=1}^{\infty}\mid\varphi g(a+x_{j})-\varphi g(a)\mid\leq\sum
_{k=1}^{\infty}\Vert\frac{1}{k!}\overset{\wedge}{d^{k}}\varphi g(a)\Vert
_{as(1;1)}\Vert(x_{j})_{j=1}^{\infty}\Vert_{w,1}^{k}=
\end{equation*}
\begin{equation*}
=\Vert(x_{j})_{j=1}^{\infty}\Vert_{w,1}\sum_{k=1}^{\infty}\Vert\frac{1}{k!}%
\overset{\wedge}{d^{k}}\varphi
g(a)\Vert_{as(1;1)}\Vert(x_{j})_{j=1}^{\infty}\Vert_{w,1}^{k-1}\leq
\end{equation*}
\begin{equation*}
\leq\Vert(x_{j})_{j=1}^{\infty}\Vert_{w,1}\sum_{k=1}^{\infty}\frac{%
e^{k}Cc^{k}\Vert\varphi\Vert}{(2ec)^{k-1}}\leq
D\Vert(x_{j})_{j=1}^{\infty}\Vert_{w,1}
\end{equation*}
for every $\varphi\in B_{F}%
%TCIMACRO{\UNICODE{0xb4}}%
%BeginExpansion
{\acute{}}%
%EndExpansion
.$ Hence
\begin{equation*}
\Vert(g(a+x_{j})-g(a))_{j=1}^{\infty}\Vert_{w,1}\leq
D\Vert(x_{j})_{j=1}^{\infty}\Vert_{w,1}\text{ }
\end{equation*}
whenever $\Vert(x_{j})_{j=1}^{\infty}\Vert_{w,p}<\delta.\blacksquare$

\begin{proposition}
\label{Geral}If $F$ has cotype $q$ and $g:E\rightarrow F$ is analytic at $%
a\in E$, then $g$ is $(q;1)$-summing at $a.$
\end{proposition}

Proof. Let $a\in E.$ Since $g$ is analytic at $a$, there exists $\delta$
such that
\begin{equation*}
\Vert(g(a+x_{j})-g(a))_{j=1}^{\infty}\Vert_{w,1}\leq
D\Vert(x_{j})_{j=1}^{\infty}\Vert_{w,1}.
\end{equation*}
Let $(x_{j})_{j=1}^{\infty}\in l_{1}^{u}(E)$ and let \ $j_{0}\in\mathbb{N}$
be such that $\Vert(x_{j})_{j=j_{0}}^{\infty}\Vert_{w,1}<\delta.$ Then
\begin{equation*}
(\sum\limits_{j=j_{0}}^{\infty}\parallel
g(a+x_{j})-g(a)\parallel^{q})^{1/q}\leq
C_{q}(F)\parallel(g(a+x_{j})-g(a))_{j=j_{0}}^{\infty}\parallel_{w,1}
\end{equation*}
\begin{equation*}
\leq D\Vert(x_{j})_{j=j_{0}}^{\infty}\Vert_{w,1}.
\end{equation*}
Obviously,
\begin{equation*}
(\sum\limits_{j=1}^{j_{0}-1}\parallel
g(a+x_{j})-g(a)\parallel^{q})^{1/q}<\infty.
\end{equation*}
Hence
\begin{equation*}
(\sum\limits_{j=1}^{\infty}\parallel
g(a+x_{j})-g(a)\parallel^{q})^{1/q}<\infty
\end{equation*}
whenever $(x_{j})_{j=1}^{\infty}\in l_{1}^{u}(E).\blacksquare$

In the real case, a slight variation of the Proposition \ref{Geral} can be
made as we see below:

\begin{proposition}
\label{ck} Let $f:E\rightarrow F$ be an application of class $C^{k}$ at $%
a\in E.$ If $\cot F\leq q$ and $\cot$ $E\leq kq$, then $f$ is$\ (q;1)$%
-summing at $a.$
\end{proposition}

Proof. Recall that if $f$ is an application of class $C^{k}$ at $a$, by
Taylor's formula there exists $B_{\delta}(a)$ such that
\begin{equation*}
\Vert f(a+x)-f(a)\Vert\leq\Vert df(a)(x)+\frac{\overset{\wedge}{d^{2}}f(a)}{%
2!}(x)+...+\frac{\overset{\wedge}{d^{k}}f(a)}{k!}(x)\Vert+\Vert x\Vert^{k}%
\text{ }\forall x\in B_{\delta}(a).
\end{equation*}
It is clear that we can consider $(x_{j})_{j=1}^{\infty}\in l_{1}^{u}(E)$ so
that \ $x_{j}\in B_{\delta}(a)$ for every $j$. Then,
\begin{equation*}
\underset{j=1}{(\overset{m}{\sum}}\Vert f(a+x_{j})-f(a)\Vert^{q})^{1/q}\leq
\end{equation*}
\begin{equation*}
\leq\lbrack\sum\limits_{j=1}^{m}(\Vert df(a)(x_{j})+\frac{\overset{\wedge }{%
d^{2}}f(a)}{2!}(x_{j})+...+\frac{\overset{\wedge}{d^{k}}f(a)}{k!}%
(x_{j})\Vert+\Vert x_{j}\Vert^{k})^{q}]^{1/q}.
\end{equation*}
Thus
\begin{equation*}
(\sum\limits_{j=1}^{m}\Vert f(a+x_{j})-f(a)\Vert^{q})^{1/q}\leq
\end{equation*}
\begin{equation*}
\leq\lbrack\sum\limits_{j=1}^{m}\Vert df(a)(x_{j})+\frac{\overset{\wedge }{%
d^{2}}f(a)}{2!}(x_{j})+...+\frac{\overset{\wedge}{d^{k}}f(a)}{k!}%
(x_{j})\Vert^{q}]^{1/q}+[\sum\limits_{j=1}^{m}(\Vert
x_{j}\Vert^{k})^{q}]^{1/q}.
\end{equation*}
Since $\cot E\leq kq$ and since $df(a),...,\overset{\wedge}{d^{k}}f(a)$ are $%
\left( q;1\right) $-summing, the proof is done.$\blacksquare$

It is not difficult to achieve the following result:

\begin{theorem}
\label{dela}If $F$ has cotype $q$, $X$ is an $\mathcal{L}_{\infty,\lambda}$
space and $f:X\rightarrow F$ is analytic at $a$, then $f$ is absolutely $%
(q;2)$-summing at $a$.
\end{theorem}

Proof.

There are $C\geq0$ and $c>0$ such that\
\begin{equation*}
\Vert\frac{1}{k!}\overset{\wedge}{d^{k}}f(a)\Vert\leq Cc^{k}\text{ for every
}k.
\end{equation*}
Thus, by (\ref{referenciar}) in Theorem \ref{77} we have
\begin{equation}
\Vert\frac{1}{k!}\overset{\wedge}{d^{k}}f(a)\Vert_{as(q;2)}\leq
C_{q}(F)K_{G}3^{\frac{k-2}{2}}\Vert\frac{1}{k!}\overset{\wedge}{d^{k}}%
f(a)\Vert\lambda^{k}\leq\text{ }C_{q}(F)K_{G}3^{\frac{k-2}{2}%
}Cc^{k}\lambda^{k}.   \label{ww1}
\end{equation}
for every $k\geq2$. \newline
For $k=1$ we know that every linear mapping from $X$ into $F$ is $(q;2)$%
-summing and it is enough, in addiction with (\ref{ww1}), to obtain positive
$C_{1}$ and $c_{1}$ so that
\begin{equation*}
\Vert\frac{1}{k!}\overset{\wedge}{d^{k}}f(a)\Vert_{as(q;2)}\leq
C_{1}c_{1}^{k}\text{ for every }k\geq1.
\end{equation*}
If $\delta_{a}$ is the radius of convergence of $f$ around $a$, then,
whenever $(x_{j})_{j=1}^{\infty}$ is such that $\Vert(x_{j})_{j=1}^{\infty}%
\Vert _{w,1}\leq\min\{\frac{1}{2c_{1}},\delta_{a}\},$ we have
\begin{equation*}
(\sum_{j=1}^{\infty}\Vert f(a+x_{j})-f(a)\Vert^{q})^{\frac{1}{q}}=\sum
_{j=1}^{\infty}(\Vert\sum_{k=1}^{\infty}\frac{1}{k!}\overset{\wedge}{d^{k}}%
f(a)(x_{j})\Vert^{q})^{\frac{1}{q}}\leq
\end{equation*}
\begin{equation*}
\leq\sum_{k=1}^{\infty}[\sum_{j=1}^{\infty}\Vert\frac{1}{k!}\overset{\wedge
}{d^{k}}f(a)(x_{j})\Vert^{q}]^{\frac{1}{q}}\leq
\end{equation*}
\begin{equation*}
\leq\sum_{k=1}^{\infty}\Vert\frac{1}{k!}\overset{\wedge}{d^{k}}f(a)\Vert
_{as(q;2)}\Vert(x_{j})_{j=1}^{\infty}\Vert_{w,2}^{k}=
\end{equation*}
\begin{equation*}
=\Vert(x_{j})_{j=1}^{\infty}\Vert_{w,2}\sum_{k=1}^{\infty}\Vert\frac{1}{k!}%
\overset{\wedge}{d^{k}}f(a)\Vert_{as(q;2)}\Vert(x_{j})_{j=1}^{\infty}%
\Vert_{w,2}^{k-1}\leq
\end{equation*}
\begin{equation*}
\leq C_{1}\Vert(x_{j})_{j=1}^{\infty}\Vert_{w,2}\sum_{k=1}^{\infty}\frac
{c_{1}^{k}}{2^{k-1}c_{1}^{k-1}}=C_{2}\Vert(x_{j})_{j=1}^{\infty}\Vert_{w,2}.
\end{equation*}

\begin{remark}
Theorem \ref{Perezz} could induce us to postulate that every mapping,
analytic at $a,$ from an $\mathcal{L}_{\infty}$ space into $\mathbb{K}$
would be $(1;2)$-summing. However, it is not true since the only absolutely $%
(1;2)$-summing linear mapping is the trivial mapping.
\end{remark}

The next example follows the same line of thought of Lemma \ref{LLema} and
Theorem \ref{dela}:

\begin{example}
If $X$ is an $\mathcal{L}_{\infty}$ space and $f:X\rightarrow\mathbb{K}$ is
a mapping, analytic at $a,$ so that $df(a)=0,$ then $f$ is $(1;2)$-summing
at $a$.
\end{example}

The reader must note that the same reasoning of Theorem \ref{dela} lead us
to the following useful Theorem:

\begin{theorem}
If the mapping $f:E\rightarrow F$ is analytic at $a\in E$ and there are $C>0$
and $c>0$ such that for each natural $n,$ $\ $%
\begin{equation}
\overset{\wedge}{d^{n}}f(a)\in\mathcal{P}_{as(s;r)}(^{n}E;F)
\end{equation}
and
\begin{equation}
\Vert\frac{1}{n!}\overset{\wedge}{d^{n}}f(a)\Vert_{_{as(s;r)}}\leq Cc^{n},%
\text{ }
\end{equation}
then $f$ is $(s;r)$-summing at $a.$
\end{theorem}

\begin{remark}
For entire holomorphic mappings we have a completer result, due to Matos
\cite{Matos2}.
\end{remark}

\begin{acknowledgement}
The author thanks Professor M.C. Matos for the suggestions, G. Botelho for
important contacts and D. P\'{e}rez-Garc\'{i}a who kindly sent him his
dissertation.
\end{acknowledgement}

\end{document}